\documentclass{article}
\usepackage{amsfonts}
\newtheorem{thm}{Theorem}
\newtheorem{prop}[thm]{Proposition}

\newtheorem{cor}[thm]{Corollary}

\newenvironment{proof}{\medskip 
\noindent {\bf Proof.}}{\hfill \rule{.5em}{1em}\mbox{}\bigskip}
\def\bea{\begin{eqnarray*}}
\def\eea{\end{eqnarray*}}
\def\bel{\begin{eqnarray}}
\def\eel{\end{eqnarray}}
\def\be{\begin{equation}}
\def\ee{\end{equation}}

\def\Bbb{\mathbb }
\def\bcp{{\mathbb C \mathbb P}}
\def\bea{\begin{eqnarray*}}
\def\eea{\end{eqnarray*}}
\newtheorem{main}{Theorem}

\title{Topology versus Chern Numbers\\ for Complex \(3\)-Folds}
\author{Claude LeBrun\thanks{Supported 
in part by  NSF grant DMS-9505744.}\\
SUNY Stony Brook} 
\date{January 28, 1998}
\begin{document}
 \maketitle
 
 \begin{abstract}
 We show by example 
  that the Chern numbers  \(  \mathbf{ 	c_{1}^{3}} \)
 and \( 	\mathbf{c_{1}c_{2}} \)
 of a complex \(3\)-fold are not determined by the 
 topology of the underlying  smooth 
 compact \(6\)-manifold.
 In fact, we observe that 
 infinitely many different values of 
 a Chern number can be achieved by 
 (integrable) complex structures on 
 a fixed  \(6\)-manifold. 
 
 \end{abstract}
 
 \section{Introduction}
 
 Suppose that \(X\) is a smooth compact oriented  \(6\)-manifold. 
 Recall that an \emph{almost-complex structure}
 on \( X \)
 means an endomorphism \( J:TX\to TX \)
 of the tangent bundle of \( X \)
with \( J^{2}=-1 \) which determines  the 
   given orientation of \(X\). 
 Such a structure makes \(TX\) into a complex
 vector bundle, so that one can speak of
 the Chern classes \( c_{j}\in H^{2j}(X,{\mathbb Z}) \)
 of \( (X,J) \), and therefore of the 
 \emph{Chern numbers}
  \begin{eqnarray*}
  \mathbf{ 	c_{1}^{3}} & = & \int_{X} c_{1}^{3}  \\
  	\mathbf{c_{1}c_{2}} & = &\int_{X} c_{1}c_{2}  \\
  	\mathbf{c_{3}} & = & \int_{X} c_{3}
  \end{eqnarray*}
  of the almost-complex manifold \( (X, J) \). 
 The only obstruction \cite{wall} to the 
 existence of an almost-complex structure 
 \(J\) on \(X\) is that \(X\) be spin\(^c\).
 This happens precisely  when the second Stiefel-Whitney 
 class \(w_{2}(X)\in H^{2}(X, {\mathbb Z}_{2})\)
 can be written as the mod-2 reduction of an element 
 of \(H^{2}(X, {\mathbb Z} )\), in which case each preimage  
 of \(w_{2}\) in 
 \(H^{2}(X, {\mathbb Z} )\)
 can be 
 realized as \(c_{1}\) for some almost-complex structure  \(J\). 
 It follows that the Chern numbers 
 \(  \mathbf{ 	c_{1}^{3}} \)
 and \( 	\mathbf{c_{1}c_{2}} \)
 of the almost-complex \((X,J)\) are certainly \emph{not} 
 topological invariants of the 
 \(6\)-manifold \(X\). For example, if 
 \( X={\mathbb C \mathbb P}_{3}, \)
 every integer  of  the form \(8j\) can be 
 realized as  \( 	\mathbf{c_{1}c_{2}} \), and every integer
 of the form \(8j^{3}\) can be
 realized as  \(  \mathbf{ 	c_{1}^{3}} \)
 for some almost-complex structure \(J\) on 
 \( {\mathbb C \mathbb P}_{3} \). 
 On the other hand, \(  c_{3} \)
 is  the Euler class of \(TX\),
 so  that \(	{\mathbf{c_{3}}}=\chi ( X )\) is actually  
 a homotopy invariant of \(X\). 
 
 In this note, we will observe that the above 
 situation  persists even if one demands 
 that the almost-complex structures under consideration
 be integrable. 
Recall that an  
 almost-complex structure \( J \) on $X$ 
 is called a  \emph{complex structure} 
 if it is integrable, in the sense of being
  locally
 isomorphic to the standard, constant-coefficient
 structure on \( {\mathbb R}^{6}={\mathbb C}^{3} \).
 The question of whether the Chern numbers 
 \(  \mathbf{ 	c_{1}^{3}} \)
 and \( 	\mathbf{c_{1}c_{2}} \)
 of a complex \(3\)-fold might actually be 
  topological invariants of the underlying
  \(6\)-manifold was raised, for example, in 
  an interesting  survey article by Okonek and van de Ven 
  \cite[p. 317]{ov}.
  
  Our principal results are as   follows: 
  
   \begin{main} \label{A}  
   There is a  compact simply connected  6-manifold  \(X\)
   which admits a sequence \(J_{m}\), \(m\in {\Bbb Z}^{+}\),
   of (integrable) complex structures
   with 
   \[ { \bf	c_{1}c_{2}} (X, J_{m})= 48m.\]
   Indeed, there are infinitely many homeotypes of 
   $X$ with 
   this property. 
  \end{main}
  
      \begin{main}\label{B}
     Let $(m,n)$ be any pair of integers. Then for 
     any integer $\tilde{n}\ll n$,
     there is a  complex projective $3$-fold 
     $(X,J)$ with Chern numbers $${\bf  c_{1}c_{2}}=24m, ~~
      {\bf c_{1}^{3}}=8n, $$ which admits 
    a second complex structure $\tilde{J}$ with 
     $${\bf  c_{1}c_{2}}=48m, ~~{\bf c_{1}^{3}}=8\tilde{n}.$$ 
     \end{main}

  \section{Infinitely Many Complex Structures}
  
  The fact that the Chern classes of a complex
  $3$-fold are not determined by the 
  topology of the underlying $6$-manifold
  was  observed long ago by Calabi \cite{calabi}. 
  While his examples all have vanishing Chern numbers, 
  they nonetheless contain the seeds of
  a natural class of examples which lead to  Theorem \ref{A}:
  
   \begin{thm}\label{k3}  For each
   positive integer $m$, 
   the 6-manifold \(X=K3 \times S^{2}\)
   admits a 
     complex structure \(J_{m}\)
   with 
   \[ { \bf	c_{1}c_{2}} (X, J_{m})= 48m\]
   and \(    {\bf 	c_{1}^{3}}(X, J_{m}) =0\).
  \end{thm}
  \begin{proof} Let $M$ denote the underlying
   oriented 4-manifold of the $K3$ surface, and 
  let \(g\) be any hyper-K\"ahler metric on 
  \(M\); such metrics exist by Yau's solution of the 
  Calabi conjecture \cite{yau}. Let \(Z\) be the twistor
  space \cite{ahs,pen} 
  of \((M, g)\), and let \(\varpi : Z\to {\mathbb C \mathbb P}_{1}\)
  be the holomorphic projection induced by the hyper-K\"ahler
  structure. Differentiably, \(\varpi\) is the trivial 
  fiber bundle with   fiber \(M\), so that \(Z\) may be thought
  of as \(X=M\times S^{2}\) equipped with a complex structure.
  
  Now let \(f_m: {\mathbb C \mathbb P}_{1}\to {\mathbb C \mathbb P}_{1}\) be 
  a holomorphic map of degree \(m-1\); for example,
  we may take \(f_m([u,v])=[u^{m-1},v^{m-1}]\), where
  \(m\) is any positive integer. We may then define
  a holomorphic family \(f_m^{*}\varpi\) of \(K3\)'s over  
  over \({\mathbb C \mathbb P}_{1}\) by pulling back the family
  \(\varpi\) via \(f_m\):
  \[\begin{array}{ccc}
  	f_m^{*} Z& \longrightarrow & Z  \\
 {\scriptstyle  f_m^{*}\varpi}	\downarrow 
 \hphantom{{\scriptstyle f_m^{*}\varpi}}&  &
 {\scriptstyle  \varpi} \downarrow \hphantom{{\scriptstyle \varpi}} \\
  	{\mathbb C \mathbb P}_{1} & \stackrel{f_m}{\longrightarrow} & {\mathbb C \mathbb P}_{1}.
  \end{array}\]
  In other words, \(f_m^{*}Z\) is the inverse image, via 
  \[ \varpi\times 1: Z\times {\mathbb C \mathbb P}_{1} 
  {\longrightarrow} {\mathbb C \mathbb P}_{1}
  \times {\bcp}_{1},\]
  of the graph of \(f_m\).  Since \(\varpi\) is 
  differentiably trivial, so is \(f_m^{*}\varpi\), and 
  \(f_m^{*}Z\) may therefore be viewed as
  \(X= K3\times S^2 \) equipped
  with a complex structure \(J_{m}\). 
  
  Now if \(\pi : Z\to M\) is the (non-holomorphic) 
   twistor projection, an explicit diffeomorphism \(Z\to X\)
   is given by \(\pi\times\varpi\), and 
   \(f_m^{*}Z\) is similarly trivialized by \(f_m^{*}\varpi
   \times f_m^{*}\pi\).
   Let \(L\subset Tf_{m}^{*}Z\) be the kernel of the derivative
   of the pulled-back twistor projection \(f_m^{*}\pi\). 
   Then \(L\) is \(J_{m}\) invariant, despite the fact that 
   \(\pi\) is not holomorphic, and so may be viewed as
   a complex line-bundle. Moreover, \(L\) may be identified
   with the pull-back of the (holomorphic) tangent bundle of
   \({\bcp}_{1}\) via \(f^{\ast}_m\varpi\), so that \(c_{1}^{2}(L)=
   0\), and hence \(p_{1}(L)=0\). If, on the other hand, 
   we use \(H\) to denote the kernel of the derivative
   of \(f_{m}^{*}\varpi\), then the underlying real bundle of 
   \(H\) is \( (f_{m}^{*}\pi)^{*}
    TM\), and so \(p_{1}(H)=f_{m}^{*}[p_{1}(M)]=
   -48 F\), where, by Poincar\'e duality,  \(F\) is 
   represented by a fiber \(S^{2}\) of \(f_m^{*}\pi\).
   It follows that \(p_{1}(f^{*}_{m}Z)=p_{1}(L)+p_{1}(H)=-48 F\).
   On the other hand, any \(K3\) has trivial canonical line bundle,
   and the fibers of \(\pi\) are \({\bcp}_{1}\)'s with normal bundle
   \({\cal   O}(1)\oplus {\cal   O}(1)\), so \(c_{1}(H)=c_{1}(L)\)
   when \(m=2\). For general \(m\), it
    follows that \(c_{1}(H)=(m-1)c_{1}(L)\), and hence that
    \[c_{1}=mc_{1}(L).\]
   We therefore have \(c_{1}^{2}=m^{2}c_1^{2}(L)=0\), so that 
   \[{\bf c_{1} c_{2}} = \frac{{\bf c_{1}^{3}}- {\bf c_{1}p_{1}}}{2}=
   24\int_{F}mc_{1}(L)=48m,\]
   and \({\bf c_{1}^{3}}=0\). 
  \end{proof}
  
  When \(m=1\), the above complex structure is simply an
  arbitrary  
  product complex structure on \(K3\times {\bcp}_{1}\),
  and so is of K\"ahler type; indeed, we may even arrange 
  for it to be 
  projective algebraic if we like. On the other hand, 
  the \(m=2\) complex structure is that of a twistor space, 
  and so is never of K\"ahler type by   Hitchin's classification 
  of K\"ahlerian twistor spaces  
  \cite{hit}. For large values of \(m\),
  one can prove something even stronger: \(J_{m}\)
  isn't even homotopic to a complex structure of 
  K\"ahler type. This is because the Todd genus
  \[ 1-h^{1}({\cal O})+ h^{2}({\cal O})=
   \chi ({\cal O}) = \frac{{\bf c_{1}c_{2}}}{24}=2m, \]
   so that \(h^{2}({\cal O})\) 
  will eventually exceed \(b_{2}(X)\),  in violation of the 
  Hodge decomposition. In the next section, we will 
  see that this phenomenon is actually quite typical.

  In order to show that there is more than one
  $6$-manifold for which  infinitely many
  different values of a Chern number are achieved by 
  (integrable) complex structures, we may now invoke the 
  standard process of {\em blowing up}. The following
facts about blow-ups $3$-folds are left as exercises for 
the reader. 
 
      \begin{prop}\label{blowup} 
     Let $(X,J)$ be any compact  complex $3$-fold, and
     let $(\hat{X}, \hat{J})$ be obtained from $(X,J)$ 
     by blowing up a point. Then $\hat{X}$ is diffeomorphic
     to the connected sum $X\# {\bcp}_{3}$, and if
     $X$ is spin, so is $\hat{X}$. Moreover,
     the
     Chern numbers of the blow-up are related to
     those of the original $3$-fold by
     \begin{eqnarray*}
 {\bf c_{1}^{3}} (\hat{X},\hat{J})& = & 
 {\bf c_{1}^{3}} (X,J) + 8\\
{\bf c_{1}c_{2}} (\hat{X},\hat{J})& = & 
   {\bf c_{1}c_{2}} (X,J) \\
  {\bf c_{3}}  (\hat{X},\hat{J})& = &  {\bf c_{3}}(X,J) +2
     \end{eqnarray*}
     \end{prop}
 
Iterated blow-ups of the previous examples thus prove the following
precise form of Thorem \ref{A}:   
          
     \begin{cor}
     For each  integer
     $n\geq 0$, the $6$-dimensional spin manifold
     $$X= (K3 \times S^{2})\# m {\bcp}_{3}$$ admits 
     a   sequence $J_{m}$ of   complex  structures   with 
     \begin{eqnarray*}
      {\bf c_{1}c_{2}}(X,J_{m}) & = & 48m  \\
      {\bf c_{1}^{3}}(X,J_{m}) & = & 8n.
     \end{eqnarray*}
     \end{cor}

  \section{K\"ahler Type}
  
  We saw in Theorem \ref{k3} that it is possible to 
  find \(6\)-manifolds with   sequences of complex structures
  for which a Chern number takes on infinitely many
  different values. On the other hand, if one
  requires the complex structures in question to 
  be of 
  {\em K\"ahler type}, one  arrives at  essentially the 
  opposite conclusion. This is 
  illustrated by our next result.
  
  \begin{thm}
  Let \(X\) be the underlying compact oriented 6-manifold of
  any K\"ahlerian $3$-fold.  
  Then  there exist infinitely many homotopy classes
  of almost-complex structures on \(X\) which {\em cannot}
  be represented by  complex structures of K\"ahler type.
  \end{thm}
  \begin{proof} By assumption, there is a K\"ahler class 
  \([\omega ] \in H^{2}(X, {\mathbb R})\) with 
  \(\ [\omega ]^{3}\neq 0\).
  Now \(H^{2}(X, {\mathbb Q})\subset 
  H^{2}(X, {\mathbb R})\) is dense, and the cup form
  is continuous, so approximating $[\omega ]$ with 
  rational classes will produce  
   classes $\alpha_{0}\in H^{2}(X, {\mathbb Q})$
  with $\alpha_{0}^{3}> 0$.  Multiplying 
  by a suitable positive integer $k$  to clear denominators, 
  we may thus obtain  
  a class $k\alpha_{0}$ which is  the image
  of an integer class \({\alpha}\in H^{2}(X, {\mathbb Z})\)
  in rational cohomology. 
  Now let \(\beta \in H^{2}(X, {\mathbb Z})\) be 
  the first Chern class of the given complex structure
  on $X$. Then \(2n \alpha + \beta\) is 
  an integer lift of \(w_{2}\), and so \cite{wall} can be 
  realized as \(c_1\) for some homotopy class \([J_n]\)
  of almost-complex structures. Now if \(J_{n}\in
  [J_{n}]\) is integrable, the Todd genus
  of \((X, J_{n})\) is 
  \[\chi ({\cal O}) = \frac{ c_{1}\cdot c_{2} }{24}=
  \frac{(2n\alpha + \beta )^{3}-(2n\alpha + \beta )\cdot p_{1}}{48}\]
  which is cubic in \(n\), with the coefficient 
  of \(n^{3}\) non-vanishing.  
  It therefore follows that there is an integer  \(N\) such that
  \(|\sum_{k}(-1)^{k}h^{0,k}| > \sum_{j}b_{j}(X)\) whenever \(|n| > N\). 
  For \(n\) in this range, the Hodge theorem must therefore fail,
  and so an integrable \(J_{n}\) could not possibly be of 
  K\"ahler type. 
  \end{proof}
  
  The reader should note that we have only used two  mild
  consequences of the K\"ahler condition: the degeneration of
  the Fr\"ohlicher spectral sequence, and the non-triviality
  of the cup form on \(H^{2}\). 
  The same argument would thus apply if one instead 
  wished to consider, 
  say, complex structures of Moishezon type.

  \section{Independence of  Chern Numbers}
  
  So far, we have seen that the Chern number ${\bf c_{1}c_{2}}$
  of a complex $3$-fold 
  is not an invariant of the underlying $6$-manifold.
  We will now see see that the 
  same is true of  ${\bf c_{1}^{3}}$.

  To this end, 
  let \(N\) be any smooth, compact oriented
  4-manifold. By \cite{taubes}, 
  the connected sum \(M= N\# k\overline{\bcp}_{2}\)
  admits anti-self-dual metrics $g$ provided that \(k\)
  is sufficiently large. The twistor space  of
   such an anti-self-dual metric is a complex 
   3-fold \((Z, J_{2})\),  the underlying $6$-manifold
   $Z$ of which is formally the 
    the fiber-wise projectivization
    ${\Bbb P}({\Bbb S}_{+})$ of the bundle of
    positive spinors on $M$. This description may 
    seem a bit paradoxical, insofar as we are 
    primarily interested in choices of $M$ which
    definitely are not spin, but it may be 
    made quite concrete by choosing a  
    spin$^{c}$ structure on $M$. This then 
    gives rise to 
    a well defined ``twisted spinor'' bundle $V_{+}$
    which formally satisfies 
     $$V_{+}={\Bbb S}_{+}\otimes L^{1/2}$$
    for  a 
    line bundle $L$ with $c_{1}(L)\cong w_{2}(M)\bmod 2$. 
    This done, we then have a canonical identification of
    $Z$ with the total space of the ${\bcp}_{1}$-bundle 
    ${\Bbb P}(V_{+})$. The naturally defined 
    complex structure $J_{2}$  then  makes 
    each fiber of  the projection $\pi : Z\to M$
    into a holomorphically embedded ${\bcp}_{1}$
    with normal bundle ${\cal O}(1)\oplus {\cal O}(1)$. 
    For more details, see \cite{ahs,hit,pen}.

   Now let us now specialize to the case in which 
    $N$ is a complex surface, and notice that
    \(M= N\# k\overline{\bcp}_{2}\) may then be thought of 
    as an iterated blow-up of $N$, and so, in particular, 
    carries a complex structure. This complex structure
    induces a spin$^{c}$ structure on $M$
    such that, for any metric $g$, the associated 
    twisted spin bundle $V_{+}$ is smoothly
    bundle-isomorphic
    to the holomorphic vector ${\cal O}\oplus K^{-1}$,
    where $K$ is the canonical line bundle of $M$.   
    Indeed, for a  Hermitian metric on $M$, there 
     is even a {\em canonical} isomorphism
     $V_{+}\cong {\cal O} \oplus K^{-1}$; 
      and, up to abstract the bundle equivalence,
     the twisted spinor bundle $V_{+}$
     is metric-independent once a spin$^{c}$ structure 
     is specified. In this way, the twistor space
     $Z$ of an anti-self-dual metric $g$ on
     $M$ is diffeomorphic to the complex manifold
     ${\Bbb P}({\cal O}\oplus K^{-1})$, 
     and so carries a second complex structure
     $J_{1}$. Notice that we do not need to assume
     any compatibility between the metric $g$ and 
     the complex structure of $M$. Also notice that
     $(Z, J_{1})$ is  projective
     algebraic (respectively,   K\"ahlerian) if $M$ is. 
     
     Let us now calculate the Chern numbers of 
     $(Z, J_{1})$. To do this, first notice that
      \(Z={\Bbb  P}({\cal   O}\oplus   K^{-1})\)
       carries two canonical
   hypersurfaces, \(\Sigma\) and \(\overline{\Sigma}\), corresponding
   to the factors of the direct sum 
  \( {\cal   O}\oplus     K^{-1} \). These are both 
   copies of the complex surface \(M\), but 
  their normal bundles are respectively \(K^{-1}\) and \(K\). 
 Moreover, the divisor \(\Sigma +\overline{\Sigma}\) precisely
 represents the vertical line bundle $L$ of $Z$. 
 We thus have
 $$c_{1}(Z, J_{1})=  \Sigma +\overline{\Sigma}+ {\pi}^{*} c_{1}(M)
 =2\Sigma . $$
 Hence
 $$  {\bf c_{1}^{3}}(Z, J_{1}) = (2\Sigma)^{3}= 8(2\chi + 3\tau )$$
     and 
     $${\bf c_{1}c_{2}}(Z, J_{1})= 2 ({\bf c_{1}^{2}}+
     {\bf c_{2}})(M) = 6 (\chi + \tau ) , $$
  where   \(\chi\) and 
   \(\tau\) are 
      the Euler characteristic and signature of \(M\),
      respectively.

     On the other hand, there is a fiber-wise antipodal map 
   which   acts anti-holomorphically  
      on $(Z, J_{2})$, so $c_{1}(J_{2})$ is Poincar\'e 
      dual to an element of $H_{4}(Z)$ which is invariant
      under this antipodal map.  We also know that 
      the integral of  $c_{1}(J_{2})$  on a fiber is 4. 
      It follows that 
      $$c_{1}(Z, J_{2})= 2\Sigma + 2\overline{\Sigma}.$$
      Since the restrictions of $J_{1}$ and $J_{2}$ 
      to a tubular neighborhood of $\Sigma$ are
      also  homotopic, we therefore deduce the
      formul{\ae} 
      \begin{eqnarray*}
       {\bf    c_{1}^{3}}  & = & 16 (2\chi + 3\tau )  \\
     {\bf  c_{1}c_{2}} & = & 12(\chi + \tau )  
   \end{eqnarray*}
     derived (with 
    opposite orientation conventions) by Hitchin  \cite{hit} in 
    greater generality. 
     For us, the point is that the invariants
     $\bf c_{1}^{3}$ and  $\bf c_{1}c_{2}$ 
     of $(Z, J_{2})$ are precisely double the corresponding 
     Chern numbers
     of $(Z, J_{1})$.

     We are now in a position to prove 
     a more precise version of Theorem \ref{B}.

     \begin{thm}
     Let $(m,n)$ be any pair of integers. Then for 
     any integer $\tilde{n}\ll n$,
     there is a spin, complex projective $3$-fold 
     $(X,J)$ with Chern numbers $${\bf  c_{1}c_{2}}=24m, ~~
      {\bf c_{1}^{3}}=8n, $$ which admits 
    a second complex structure $\tilde{J}$ with 
     $${\bf  c_{1}c_{2}}=48m, ~~{\bf c_{1}^{3}}=8\tilde{n}.$$ 
      If $m > 0$, moreover, 
      we may even arrange for $X$ to be 
     simply connected.
     \end{thm}
     
 \begin{proof}
 For each integer $m$, we begin by choosing a complex
 surface $N(m)$ with   Todd genus $(\chi + \tau )/4=m$. 
 For example, if $m\leq 1$, let us take $N$ to be 
 $C\times {\bcp}_{1}$, where $C$ is a Riemann
 surface of genus $1-m$. On the other hand, if 
 $m > 1$, we may let $N$ be the minimal resolution
 of $(E\times C)/{\Bbb Z}_{2}$, where $E$ is an
 elliptic curve, $C$ is a hyperelliptic curve
 of genus $m-1$, and ${\Bbb Z}_{2}$ acts simultaneously on both
 factors by the Weierstrass involution. Notice that 
 our choice of 
 $N(m)$ is simply connected when $m > 0$,
 and that, incidentally,  this is the best
 one can do in principal.  
 
 Now, for each $m$, let $k_0(m)$ be chosen so that 
$N(m)\# k\overline{\bcp}_{2}$ admits anti-self-dual metrics
for each $k \geq k_0(m)$. By Taubes' theorem
 \cite{taubes}, such an integer $k_0(m)$ exists. 
Moreover, with the above choices, we may even take $k_0(m)= 0$
for $m < 0$, $k_0(0)=6$, $k_0(1)=14$, and $k_0(2)= 3$
\cite{L2,LS2,DF}. 

Now let $(m,n)$ be any pair of integers, and let 
$\tilde{n}$ be any integer such that
$$\tilde{n} \leq \min (n-k_{0}(m)+ {\bf c_{1}^{2}}(N(m)), 2n).$$
We may then define    integers $k\geq k_{0}(m)$ and
$\ell \geq 0$ by 
\begin{eqnarray*}
	k & = & n-\tilde{n}+ {\bf c_{1}^{2}}(N(m)),   \\
	\ell & = & 2n-\tilde{n}.
\end{eqnarray*}
Let $Z(k,m)$ be the twistor space
of an anti-self-dual metric on 
$M=N(m)\# k\overline{\bcp}_{2}$, and  let
  $$X(k,\ell, m) = Z(k,m) \# \ell \bcp_{3}.$$
  Notice that $X$ is a spin manifold. Moreover, it
  comes equipped with two different complex structures.

   First, because $M=N(m)\# k\overline{\bcp}_{2}$ 
  is a projective algebraic surface,  $Z$ carries a 
  projective algebraic complex structure $J_{1}$, and 
  $X$ may then be identified with the blow-up of
  $(Z, J_{1})$ at $\ell$ points. Let us denote
   this   complex structure
  on $X$ by $J$.  Using Proposition \ref{blowup} and the
  above computations, we therefore have 
  \begin{eqnarray*}
       {\bf    c_{1}^{3}} (X,J) & = & 8 (2\chi + 3\tau ) (M) +8\ell= 
       8n\\
     {\bf  c_{1}c_{2}}(X,J)  & = & 6(\chi + \tau ) (M) = 24m. 
  \end{eqnarray*}
  On the other hand, each $Z$ also admits its twistor complex
  structure $J_{2}$, and we may  instead choose
  to think of $X$ as the blow-up of this twistor space at
  $\ell$ points. Let us denote the corresponding  complex structure 
  on $X$ by $\tilde{J}$. Thus 
   \begin{eqnarray*}
       {\bf    c_{1}^{3}} (X,\tilde{J})  & = & 16 (2\chi + 3\tau ) (M) +8\ell= 
       8\tilde{n},\\
     {\bf  c_{1}c_{2}}(X,\tilde{J})   & = & 12(\chi + \tau ) (M) = 48m, 
  \end{eqnarray*}
  as claimed. 
 \end{proof}

 Since the Todd genus of any complex $3$-fold is 
 given by ${\bf  c_{1}c_{2}}/24$, this result realizes all
 possible values of ${\bf  c_{1}c_{2}}$. The divisibility 
 of ${\bf  c_{1}^{3}}$ by $8$ is also a necessary for
 $X$ to be spin, so the result is also essentially optimal in this
 respect.

  On the other hand, we have chosen to ignore ${\bf c_{3}}$,
which is determined by $(m,n,\tilde{n})$ in these
examples. The abundance of rational curves 
 also forces all 
 our $3$-folds all have Kodaira dimension $-\infty$. 
And finally, most of our manifolds are in no 
sense  minimal. It would  obviously be of great interest
to produce new examples which overcome
 these limitations.

 \bigskip 
  
\noindent  {\bf Acknowledgments.}
The author would like to thank Ron Stern for 
drawing his attention to  the problem and
 pointing out 
several important references. He would also like to 
the thank the Mathematics Department of Harvard
University for its hospitality during the inception
of this work.

\end{document}